\documentclass[a4paper,12pt]{amsart}
\usepackage{amssymb,amsthm}

\begin{document}

\newcommand{\opp}{\bowtie }
\newcommand{\po}{\text {\rm pos}}
\newcommand{\supp}{\text {\rm supp}}
\newcommand{\End}{\text {\rm End}}
\newcommand{\diag}{\text {\rm diag}}
\newcommand{\Lie}{\text {\rm Lie}}
\newcommand{\Ad}{\text {\rm Ad}}
\newcommand{\car}{\mathcal R}
\newcommand{\Tr}{\rm Tr}
\newcommand{\Spec}{\text{\rm Spec}}

\newtheorem*{th1}{Proposition 2.4}
\newtheorem*{th2}{Proposition 2.6 {\rm (\cite[11.2(c)]{L1})}}
\newtheorem*{th3}{Lemma 3.2}
\newtheorem*{th4}{Corollary 3.7}
\newtheorem*{th5}{Lemma 3.9}
\newtheorem*{th6}{Corollary 3.10}
\newtheorem*{th7}{Lemma 3.11}
\newtheorem*{th8}{Lemma 3.12}
\newtheorem*{th9}{Proposition 4.1}
\newtheorem*{th10}{Corollary 4.5}
\newtheorem*{th11}{Corollary 4.6}
\newtheorem*{th12}{Corollary 4.7}
\newtheorem*{th13}{Proposition 5.2}
\newtheorem*{th14}{Lemma 5.4}
\newtheorem*{th15}{Corollary 5.5}
\newtheorem*{th16}{Proposition 5.9}
\newtheorem*{th17}{Corollary 5.12}

\theoremstyle{remark}
\newtheorem*{rmk}{Remark}

\def\ge{\geqslant}
\def\le{\leqslant}
\def\a{\alpha}
\def\b{\beta}
\def\c{\chi}
\def\g{\gamma}
\def\G{\Gamma}
\def\d{\delta}
\def\D{\Delta}
\def\L{\Lambda}
\def\e{\epsilon}
\def\et{\eta}
\def\io{\iota}
\def\o{\omega}
\def\p{\pi}
\def\ph{\phi}
\def\ps{\psi}
\def\r{\rho}
\def\s{\sigma}
\def\t{\tau}
\def\th{\theta}
\def\k{\kappa}
\def\l{\lambda}
\def\z{\zeta}
\def\v{\vartheta}
\def\x{\xi}
\def\i{^{-1}}

\def\mapright#1{\smash{\mathop{\longrightarrow}\limits^{#1}}}
\def\mapleft#1{\smash{\mathop{\longleftarrow}\limits^{#1}}}
\def\mapdown#1{\Big\downarrow\rlap{$\vcenter{\hbox{$\scriptstyle#1$}}$}}

\def\ca{\mathcal A}
\def\cb{\mathcal B}
\def\cc{\mathcal C}
\def\cd{\mathcal D}
\def\ce{\mathcal E}
\def\cf{\mathcal F}
\def\cg{\mathcal G}
\def\ch{\mathcal H}
\def\ci{\mathcal I}
\def\cj{\mathcal J}
\def\ck{\mathcal K}
\def\cl{\mathcal L}
\def\cm{\mathcal M}
\def\cn{\mathcal N}
\def\co{\mathcal O}
\def\cp{\mathcal P}
\def\cq{\mathcal Q}
\def\car{\mathcal R}
\def\cs{\mathcal S}
\def\ct{\mathcal T}
\def\cu{\mathcal U}
\def\cv{\mathcal V}
\def\cw{\mathcal W}
\def\cz{\mathcal Z}
\def\cx{\mathcal X}
\def\cy{\mathcal Y}

\def\tz{\tilde Z}
\def\tl{\tilde L}
\def\tc{\tilde C}
\def\ta{\tilde A}
\def\tx{\tilde X}

\title[The character sheaves on the group compactification]
{The character sheaves on the group compactification}
\author{Xuhua He}
\address{Department of Mathematics, Massachusetts Institute of Technology, Cambridge, MA 02139, USA}%
\thanks{The author thanks the National Science Foundation for its support through grant DMS-0243345 (principal investigator: George Lusztig).}
\email{xuhua@mit.edu}


\begin{abstract} We give a definition of character
sheaves on the group compactification which is equivalent to
Lusztig's definition in [Moscow Math. J. 4 (2004) 869-896]. We
also prove some properties of the character sheaves on the group
compactification.
\end{abstract}
\maketitle

\section*{Introduction}
Let $G$ be a connected reductive algebraic group over an
algebraically closed field. In \cite{L3}, Lusztig introduced a
class of $G \times G$-varieties. We denote the varieties by
$\tz_{J, y, \d}$. The precise definition can be found in 1.3. The
group $G$ acts diagonally on $\tz_{J, y, \d}$. Lusztig introduced
a partition of $\tz_{J, y, \d}$ into finitely many $G$-stable
pieces. The $G$-orbits on each piece are in one-to-one
correspondence with the conjugacy classes of a certain (smaller)
reductive group.

To each character sheaf on the (smaller) reductive group, we
associate a $G$-equivariant simple perverse sheaf on the
$G$-stable piece. We call it a character sheaf on the $G$-stable
piece. Its perverse extension to $\tz_{J, y, \d}$ is a
$G$-equivariant simple perverse sheaf on $\tz_{J, y, \d}$. The
simple perverse sheaves obtained in this way are called the
``parabolic character sheaves''. (This is a generalization of the
character sheaves on the group.)

The definition above is one of the two equivalent definitions in
\cite{L3}. The other one imitates Lusztig's definition of
character sheaves of the group. Roughly speaking, the second
definition help us to understand the perverse extensions of
character sheaves on the $G$-stable pieces to $\tz_{J, y, \d}$. A
consequence of the coincidence of the two definition is that the
parabolic character sheaves have the following property: any
composition factor of the perverse cohomology of the restriction
of a parabolic character sheaf to a $G$-stable piece is a
character sheaf on that piece.

The varieties $\tz_{J, y, \d}$ include more or less as a special
case the strata of the group compactification $\overline{G^1}$ of
$G$, here $G$ is adjoint. Therefore, we obtain a partition of
$\overline{G^1}$ into finitely many $G$-stable pieces. Lusztig
defined the character sheaves on $\overline{G^1}$ to be the
perverse extensions of the character sheaves on the $G$-stable
pieces of $\overline{G^1}$. One may expect that the character
sheaves on the group compactification have the property analogous
to the property of parabolic character sheaves that we mentioned
above.

To achieve this goal, we need to understand the perverse extension
to $\overline{G^1}$. The main propose of this paper is to
introduce an equivalent definition that helps us to understand it.

Before introducing the definition, we first recall the definition
of character sheaves of the group $G$. Denote by $B$ the Borel
subgroup of $G$. The group $G$ is decomposed into finitely many $B
\times B$-orbits. Starting with certain local systems on a $B
\times B$-orbits, we obtain some simple perverse sheaves on $G$.
(For more details, see \cite[no. 5]{MS1}.) Although this is not
Lusztig's original definition, one can see that they are
equivalent directly from the definition.

Now let us come back to the group compactification
$\overline{G^1}$. There are finitely many $B \times B$-orbits on
$\overline{G^1}$. The closure relation of the $B \times B$-orbits
and the local systems on the orbits were studied by Springer in
\cite{S1}. These local systems play the same role in our
definition of character sheaves on $\overline{G^1}$ as the local
systems on the $B \times B$-orbits of $G$ in the definition of
character sheaves on $G$. Then we show that the simple perverse
sheaves on $\overline{G^1}$ obtained in this way coincide with
those obtained from Lusztig's definition.

Now let us discuss about the ideas of the proof of the coincidence
of our definition and Lusztig's definition. We denote by
$\cd(\overline{G^1})$ the derived category of constructible
sheaves on $\overline{G^1}$. To each element in
$\cd(\overline{G^1})$, we associate its support, which is a closed
subvariety. The counterpart (on the ``level of varieties'') of the
equivalence of these definitions (on the ``level of sheaves'') is
the relation between the (closure of) $G$-stable pieces and the
(closure of) $B \times B$-orbits. The relation was discussed in
\cite{H1} and \cite{H2}. In \cite[2.7]{H1}, we showed that each
$G$-stable piece is the minimal $G$-stable subvariety that
contains a particular $B \times B$-orbit. In \cite[4.3]{H2}, we
gave an inductive way to determine in which $G$-stable piece an
element of a $B \times B$-orbit is contained. In this paper, we
``lift'' these results from the ``level of varieties'' to the
``level of sheaves''.

The content of this paper is arranged as follows. In section 1, we
recall the definitions of $\tz_{J, y, \d}$ and the $G$-stable
pieces on it. We also discuss some properties of the $G$-stable
pieces. In section 2, we recall the definition of the character
sheaves on the group and discuss some properties of them. In
section 3, we prove our key lemma. In section 4, we first
introduce a new definition of the character sheaves of group
compactification and prove that the character sheaves on the group
compactification have the ``nice'' property according to our
definition. As a consequence of the ``nice'' property, we show
that our definition is equivalent to Lusztig's. We also obtain a
property about the central characters. Our approach can also be
generalized to parabolic character sheaves. In section 5, we
discuss some results on the parabolic character sheaves. We also
obtain a new proof of the coincidence of Lusztig's two definitions
of parabolic character sheaves.

We thank Lusztig for many useful discussions on character sheaves.
We also thank Springer for some advice and comments on an earlier
version of this paper.

\section{The $G$-stable pieces}

\subsection{} Let $G$ be a connected,
reductive algebraic group over an algebraically closed field
$\bold k$. Let $B$ be a Borel subgroup of $G$, $T \subset B$ be a
maximal torus and $B^-$ be the opposite Borel subgroup. Let $W$ be
the corresponding Weyl group and $(s_i)_{i \in I}$ the set of
simple reflections. For $w \in W$, we denote by $\supp(w) \subset
I$ the set of simple roots whose associated simple reflections
occur in some (or equivalently, any) reduced decomposition of $w$.
We also choose a representative $\dot w$ of $w$ in $G$.

For $J \subset I$, let $W_J$ be the subgroup of $W$ generated by
$\{s_j \mid j \in J\}$ and $W^J$ (resp. $^J W$) be the set of
minimal length coset representatives of $W/W_J$ (resp. $W_J
\backslash W$). Let $w^J_0$ be the unique element of maximal
length in $W_J$. (We simply write $w_0$ for $w^I_0$.) For $J, K
\subset I$, we write $^J W^K$ for $^J W \cap W^K$.

For $J \subset I$, let $P_J \supset B$ be the standard parabolic
subgroup defined by $J$ and $P^-_J \supset B^-$ be the opposite of
$P_J$. Let $L_J=P_J \cap P^-_J$. For any $J \subset I$, let
$\cp_J$ be the set of parabolic subgroups conjugate to $P_J$. We
simply write $\cb$ for $\cp_{\varnothing}$. For $J, K \subset I$,
$w \in {}^J W^k$ and $P \in \cp_J$, $Q \in \cp_K$, we write
$\po(P, Q)=w$ if $g P g \i=P_J$, $g Q g \i=\dot w P_K \dot w \i$
for some $g \in G$.

For any parabolic subgroup $P$, we denote by $U_P$ its unipotent
radical. We simply write $U$ for $U_B$ and $U^-$ for $U_{B^-}$. We
denote by $H_P$ the inverse image of the connected center of $P /
U_P$ under $P \rightarrow P / U_P$.

For any closed subgroup $H$ of $G$, we denote by $H_{diag}$ the
image of the diagonal embedding of $H$ in $G \times G$.

\subsection{} Let $\hat{G}$ be a possibly disconnected
reductive algebraic group over $\bold k$ with identity component
$G$. For $g \in \hat{G}$ and $H \subset G$, we write ${}^g H$ for
$g H g \i$. Let $G^1$ be a connected component of $\hat{G}$. There
exists an isomorphism $\d: W \rightarrow W$ such that $\d(I)=I$
and ${}^g P \in \cp_{\d(J)}$ for $g \in G^1$ and $P \in \cp_J$. We
choose $g_0 \in G^1$ in the same way as \cite[1.2]{H2}. In
particular, if $G^1=G$, then $g_0=1$.

\subsection{} Let $J, J' \subset I$ and $y \in {}^{J'} W^{\d(J)}$
such that $y \d(J)=J'$. For $P \in \cp_J$ and $Q \in \cp_{J'}$,
set $A_y(P, Q)=\{g \in G^1 \mid \po(Q, {}^g P)=y\}$. Set $$\tz_{J,
y, \d}=\{(P, Q, \g) \mid P \in \cp_J, Q \in \cp_{J'}, \g \in
U_{P'} \backslash A_y(P, Q) /U_P\}.$$

Define the $G \times G$ action on $\tz_{J, y, \d}$ by $(g_1, g_2)
\cdot (P, Q, \g)=({}^{g_2} P, {}^{g_1} Q, g_1 \g g_2 \i)$. It is
easy to see that $G \times G$ acts transitively on $\tz_{J, y,
\d}$.

Set $$\tilde{h}_{J, y, \d}=(P_J, {}^{\dot y \i} P_{J'}, U_{^{\dot
y \i} P_{J'}} g_0 U_{P_J}) \in \tz_{J, y, \d}.$$

For $w, v \in W$, set $$\widetilde{[J, w, v]}_{y, \d}=(B \times B)
(\dot w, \dot v) \cdot \tilde{h}_{J, y, \d}.$$

It is easy to see that $\widetilde{[J, w u, v  \d \i(u)]}_{y,
\d}=\widetilde{[J, w, v]}_{y, \d}$ for $u \in W_{\d(J)}$ and
$$\tz_{J, y, \d}=\bigsqcup_{w \in W^{\d(J)}, v \in W} \widetilde{[J, w, v]}_{y,
\d}.$$

For $w \in W^{\d(J)}$, set $$\tz^w_{J, y, \d}=G_{diag} \cdot
\widetilde{[J, w, 1]}_{y, \d}.$$

By \cite[1.3 \& 1.7]{H2}, we have that
$$\tz_{J, y, \d}=\bigsqcup_{w \in W^{\d(J)}} \tz^w_{J, y, \d}.$$

The subvarieties $\tz^w_{J, y, \d}$ are called the $G$-stable
pieces of $\tz_{J, y, \d}$.

\subsection{} For $w \in W^{\d(J)}$, set $$I(J, w, \d)=\max\{K \subset J \mid
w \d(K)=K\}.$$ In the rest of this subsection, we fix $w \in
W^{\d(J)}$ and write $K$ for $I(J, w, \d)$. By \cite[1.10]{H2}, we
have the following results.

(1) The map $G \times \bigl((P_K)_{diag} (L_K \dot w, 1) \cdot
\tilde h_{J, y, \d} \bigr) \rightarrow \tz^w_{J, y, \d}$ defined
by $(g, z) \mapsto (g, g) \cdot z$ induces an isomorphism $$G
\times_{P_K} \bigl((P_K)_{diag} (L_K \dot w, 1) \cdot \tilde h_{J,
y, \d}) \bigr) \cong \tz^w_{J, y, \d},$$ where $P_K$ acts on $G$
on the right and acts on $(P_K)_{diag} (L_K \dot w, 1) \cdot
\tilde h_{J, y, \d})$ diagonally.

(2) The map $P_K \times L_K \dot w g_0 \rightarrow (P_K)_{diag}
(L_K \dot w, 1) \cdot \tilde h_{J, y, \d}$ defined by $(p, l)
\mapsto (p l g_0 \i, p) \cdot \tilde h_{J, y, \d}$ for $p \in P_K$
and $l \in L_K \dot w g_0$ induces an affine space bundle map
$$P_K \times_{L_K} L_K \dot w g_0 \rightarrow (P_K)_{diag} (L_K
\dot w, 1) \cdot \tilde h_{J, y, \d},$$ where $L_K$ acts on $G$ on
the right and acts on $L_K \dot w g_0$ by conjugation.

\

Therefore, we obtain an affine space bundle map $$\tilde \pi^w_{J,
y, \d}: G \times_{L_K} L_K \dot w g_0 \rightarrow \tz^w_{J, y,
\d}$$ which sends $(g, l)$ to $(g l g_0 \i, g) \cdot \tilde h_{J,
y, \d}$ for $g \in G$ and $l \in L_K \dot w g_0$.

Let $\v: \tz^w_{J, y, \d} \rightarrow \tz^w_{K, w, \d}$ be the
morphism defined by $$\v \bigl((g l g_0 \i, g) \cdot \tilde h_{J,
y, \d} \bigr)=(g l g_0 \i, g) \cdot \tilde h_{K, w, \d}$$ for $g
\in G$, $p \in P_K$ and $l \in L_K \dot w g_0$. By \cite[1.6 \&
1.10]{H2}, $\v$ is well defined. Moreover, we have the following
commuting diagram

$$
\begin{matrix} G \times_{L_K} L_K \dot w g_0 & \mapright {\tilde \pi^w_{J, y, \d}}
& \tz^w_{J, y, \d} \\ \Big\Vert & & \mapdown\v \\ G \times_{L_K}
L_K \dot w g_0 & \mapright {\tilde \pi^w_{K, w, \d}} & \tz^w_{K,
w, \d}.
\end{matrix}$$

Note that $w \d(K)=K$. It is easy to see that $L_K \dot w g_0$ is
a connected component of $N_{\hat{G}}(L_K)$. However, in general
$L_K \dot w g_0$ is not the identity component of
$N_{\hat{G}}(L_K)$ even if $\hat{G}=G=G^1$. We will see in 4.3
that the character sheaves on a (possibly) disconnected group are
involved when studying character sheaves on the group
compactification even if $\hat{G}=G=G^1$.

\subsection{} Let $J, J' \subset I$ and $y \in {}^{J'} W^{\d(J)}$
such that $y \d(J)=J'$. Set $$Z_{J, y, \d}=\{(P, Q, \g) \mid P \in
\cp_J, Q \in \cp_{J'}, \g \in H_{P'} \backslash A_y(P, Q) /H_P\}$$
with the $G \times G$ action defined in the same way as $\tz_{J,
y, \d}$.

As in \cite[11.19]{L3}, we may identify groups $H_P / U_P$ (with
$P \in \cp_J$) with a single torus $\D_J$ independent of the
choice of $P$. Now $\D_J$ acts freely on $\tz_{J, y, \d}$ by $t:
(P, Q, \g) \mapsto (P, Q, \g z)$ where $z \in H_P$ represents $t
\in D_J$. Then we may identify $Z_{J, y, \d}$ with $\D_J
\backslash \tz_{J, y, \d}$ as $G \times G$-varieties.

Set $h_{J, y, \d}=(P_J, {}^{\dot y \i} P_{J'}, H_{^{\dot y \i}
P_{J'}} g_0 H_{P_J}) \in Z_{J, y, \d}$.

For $w, v \in W$, set $[J, w, v]_{y, \d}=(B \times B) (\dot w,
\dot v) \cdot h_{J, y, \d}$. For $w \in W^{\d(J)}$, set $Z^w_{J,
y, \d}=G_{diag} \cdot [J, w, 1]_{y, \d}$. Then
$$Z_{J, y, \d}=\bigsqcup_{w \in W^{\d(J)}, v \in W} [J, w, v]_{y,
\d}=\bigsqcup_{w \in W^{\d(J)}} Z^w_{J, y, \d}.$$

The subvarieties $Z^w_{J, y, \d}$ are called the $G$-stable pieces
of $Z_{J, y, \d}$.

Now fix $w \in W^{\d(J)}$. Set $K=I(J, w, \d)$. Then as a
consequence of 1.4, we obtain the following commuting diagram
$$\begin{matrix} G \times_{L_K} L_K \dot w g_0 /Z^0(L_J) & \mapright{} & Z^w_{J, y, \d} \\
\Big\Vert & & \mapdown\v \\ G \times_{L_K} L_K \dot w g_0 /
Z^0(L_J) & \mapright{} & Z^w_{K, w, \d},
\end{matrix}$$ where $Z^0(L_J)$ is the connected center of $L_J$
and $\v$ is defined in the similar way as in 1.4.

\subsection{} In this subsection, we assume that $G$ is adjoint. The
compactification $\overline{G^1}$ of $G^1$ is the $G \times G$
variety which is isomorphic to the wonderful compactification of
$G$ as a variety and where the $G \times G$ action is twisted by
$G \times G \rightarrow G \times G$, $(g, g') \mapsto (g, g_0 g'
g_0 \i)$. By \cite[2.1]{H2}, we have that
\begin{align*} \overline{G^1} &=\bigsqcup_{J \subset I} Z_{J, w_0 w^{\d(J)}_0,
\d}=\bigsqcup_{J \subset I} \bigsqcup_{w \in W^{\d(J)}, v \in W}
[J, w, v]_{w_0 w^{\d(J)}_0, \d} \\ &=\bigsqcup_{J \subset I}
\bigsqcup_{w \in W^{\d(J)}} Z^w_{J, w_0 w^{\d(J)}_0, \d}.
\end{align*}

We write $h_{J, w_0 w^{\d(J)}_0, \d}$ as $h_J$, $[J, w, v]_{w_0
w^{\d(J)}_0, \d}$ as $[J, w, v]$ and $Z^w_{J, w_0 w^{\d(J)}_0,
\d}$ as $Z^w_J$. We call $Z^w_J$ the $G$-stable pieces of
$\overline{G^1}$. We denote by $\pi^w_J$ the morphism $G
\times_{L_{I(J, w, \d)}} L_{I(J, w, \d)} \dot w g_0 /Z^0(L_J)
\rightarrow Z^w_J$ in 1.5.

\section{The character sheaves on $G^1$}

\subsection{} We follow the notation of \cite{BBD}. Let $X$ be an
algebraic variety over $\bold k$ and $l$ be a fixed prime number
invertible in $\bold k$. We write $\cd(X)$ instead of $\cd^b_c(X,
\bar{\mathbb Q}_l)$. If $K \in \cd(X)$ and $A$ is a simple
perverse sheaf on $X$, we write $A \dashv K$ if $A$ is a
composition factor of ${}^p H^i(K)$ for some $i \in \mathbb Z$.
For $A, B \in \cd(X)$, we write $A=B[\cdot]$ if $A=B[m]$ for some
$m \in \mathbb Z$.

Let $K, K_1, K_2, \cdots, K_n \in \cd(X)$. We say that $K \in
<K_1, K_2, \cdots, K_n>$ if there exist $m \ge n$ and $K_{n+1},
K_{n+2}, \cdots, K_m \in \cd(X)$ such that $K_m=K$ and for each
$n+1 \le i \le m$, there exists $1 \le j, k<i$ and $n_j, n_k \in
\mathbb Z$, such that $(K_j[n_j], K_i, K_k[n_k])$ is a
distinguished triangle in $\cd(X)$. In this case, if
$H^{\cdot}_c(X, K_i)=0$ for all $i$, then $H^{\cdot}_c(X, K)=0$.
Moreover, if $A \dashv K$, then $A \dashv K_i$ for some $1 \le i
\le n$.

If $f: X \rightarrow Y$ is a morphism of algebraic varieties, we
have functors $f_!$ and $f^*$ between the derived categories
$\cd(X)$ and $\cd(Y)$. Therefore

(1) if $A, A_1, A_2, \cdots, A_n \in \cd(X)$ with $A \in <A_1,
A_2, \cdots, A_n>$, then $f_! A \in <f_! A_1, f_! A_2, \cdots, f_!
A_n>$;

(2) if $B, B_1, B_2, \cdots, B_n \in \cd(Y)$ with $B \in <B_1,
B_2, \cdots, B_n>$, then $f^* B \in <f^* B_1, f^* B_2, \cdots, f^*
B_n>$.

Let $H$ be a connected algebraic group and $X$, $Y$ be varieties
with a free $H$-action on $X \times Y$. Denote by $X \times^H Y$
the quotient space. For $K_1 \in \cd(X)$ and $K_2 \in \cd(Y)$ such
that $K_1 \boxtimes K_2$ is $H$-equivariant, we denote by $K_1
\odot K_2$ be the element in $\cd(X \times^H Y)$ whose inverse
image under $X \times Y \rightarrow X \times^H Y$ is $K_1
\boxtimes K_2$.

\subsection{} Let $p$ be the characteristic of $k$ and $\mathbb Z_{(p)}$ be the
ring of rational numbers with denominator prime to $p$ (In
particular, $\mathbb Z_{(p)}=Q$ if $p=0$). Set $\hat{X}=\mathbb
Z_{(p)} \otimes_{\mathbb Z} X/ 1 \otimes_{\mathbb Z} X$, where $X$
is the character group of $T$. If necessary, we write
$\hat{X}(T)$.

Let $\ck(T)$ be the set of isomorphism classes of Kummer local
systems on $T$, i. e., the set of isomorphism classes of
$\bar{\mathbb Q}_l$-local systems $\cl$ of rank one on $T$, such
that $\cl^{\otimes m} \cong \bar{\mathbb Q}_l$ for some integer $m
\ge 1$ invertible in $\bold k$. By \cite[2.1]{MS1}, we may
identify $\ck(T)$ with $\hat{X}$. For $\xi \in \hat{X}$, we denote
by $\cl_{\xi}$ the corresponding Kummer local system.

For $J \subset I$, set $T_J=T/Z^0(L_J)$. We identify
$\hat{X}(T_J)$ with a subgroup of $\hat{X}$, see \cite[section
1]{S1}.

\

In the rest of this section, we recall the definition and some
properties of the character sheaves on $G^1$. We follow the
approach in \cite{MS1}. (In [loc. cit.] Mars and Springer dealt
with the case where $G^1=G$. The general case can be treated in a
similar way.)

\subsection{} By the Bruhat decomposition, we have that
$G^1=\bigsqcup_{w \in W} B \dot w g_0 B$. For $\xi \in \hat{X}$,
let $\cl_{\xi, w, \d}$ be the inverse image of $\cl_{\xi}$ under
$B \dot w g_0 B \rightarrow T$, $u \dot w g_0 t u' \mapsto t$ for
$u, u' \in U$ and $t \in T$. Then $\cl_{\xi, w, \d}$ is a tame
local system on $B \dot w g_0 B$. We denote by $A_{\xi, w, \d}$
the perverse extension of $\cl_{\xi, w, \d}$ to $G^1$. If $w
\d(\xi)=\xi$, then $\cl_{\xi, w, \d}$ and $A_{\xi, w, \d}$ are
equivariant for the conjugation action of $B$.

Define the $B$-action on $G \times G^1$ by $b(g, h)=(g b \i, b h
b\i)$. Let $G \times_B G^1$ be the quotient space. The map $G
\times G^1 \rightarrow G^1 $, $(g, h) \mapsto g h g \i$ induces a
proper morphism $\g: G \times_B G^1 \rightarrow G^1$. By standard
arguments, we have the following result.

\begin{th1} Let $\xi \in \hat{X}$ and $A$ be a simple perverse
sheaf on $G^1$. The following conditions on $A$ are equivalent:

(i) $A \dashv (\g)_! (\bar{\mathbb Q}_l \odot A_{\xi, w, \d})$ for
some $w \in W$ with $w \d(\xi)=\xi$.

(ii) $A \dashv (\g)_! (\bar{\mathbb Q}_l \odot \cl_{\xi, w, \d})$
for some $w \in W$ with $w \d(\xi)=\xi$.
\end{th1}

\subsection*{2.5} Let $\cc_{\xi}(G^1)$ be the set of
(isomorphism classes) of simple perverse sheaves on $G^1$ which
satisfies the equivalent conditions (i)-(ii) with respect to
$\xi$. The simple perverse sheaves on $G^1$ which belong to
$\cc_{\xi}(G^1)$ for some $\xi \in \hat{X}$ are called character
sheaves on $G^1$; they (or their isomorphism classes) form a set
$\cc(G^1)$.

\begin{th2} Let $\xi, \eta \in \hat{X}$ with $\eta \notin W \xi$. Then
$$\cc_{\xi}(G^1) \cap \cc_{\eta}(G^1)=\varnothing.$$
\end{th2}

\section{The key lemma}

In this section and the next section, we assume that $G$ is
adjoint.

\subsection{} Let $J \subset I$ and $w, v \in W$. Let $p_{J, w, v}: B
\times B \rightarrow [J, w, v]$ be the morphism defined by $(b_1,
b_2) \mapsto (b_1 \dot w, b_2 \dot v) \cdot h_J$ for $b_1, b_2 \in
B$. By 1.3, $p_{J, w u, v \d \i(u)}=p_{J, w, v}$ for $u \in
W_{\d(J)}$. The morphism $B \times B \rightarrow T_J$, $(t_1 u_1,
t_2 u_2) \mapsto (\dot w g_0) \i t_1 (\dot w g_0) \dot v \i t_2 \i
\dot v Z^0(L_J)$ factors through a morphism $pr_{J, w, v}: [J, w,
v] \rightarrow T_J$.

For $\xi \in \hat{X}(T_J) \subset \hat{X}$, $pr_{J, w, v}^*
\cl_{\xi}=\cl_{\xi, J, w, v}$ is a (tame) local system on $[J, w,
v]$. By \cite[4.1.2]{MS1} and \cite[2.2.3]{MS2}, $\cl_{\xi, J, w,
v}$ has weight $(w \d(\xi), -v \xi)$ for the $B \times B$-action
and the rank one (tame) local system on $[J, w, v]$ which have a
weight for the $B \times B$-action are of the form $\cl_{\xi, J,
w, v}$ for some $\xi \in \hat{X}(T_J)$. Moreover, $\cl_{\xi, J, w
u, v \d \i(u)}=\cl_{\d \i(u) \xi, J, w, v}$ for $u \in W_{\d(J)}$.
If $w \d(\xi)=v \xi$, then $\cl_{\xi, J, w, v}$ is equivariant for
the diagonal action of $B$.

We denote by $A_{\xi, J, w, v}$ the perverse extension of
$\cl_{\xi, J, w, v}$ to $\overline{G^1}$. By \cite[2.2]{MS2}, for
$i \in \mathbb Z$, $\ch^i A_{\xi, J, w, v} \mid_{[J', w', v']}$ is
a direct sum of local system $\cl_{\eta, J', w', v'}$, where $\eta
\in \hat{X}(T_{J'})$ and $w \d(\xi)=w' \d(\eta)$, $v \xi=v' \eta$.
It is zero unless $[J', w', v']$ is contained in the closure of
$[J, w, v]$.

In other words, we have the following result.

\begin{th3} Let $J \subset I$, $w, v \in W$ and $\xi \in \hat{X}(T_J)$.
Set \begin{align*} \ci_{J, w, v, \xi}=\Bigl\{(J', w', v', \eta)
\mid & J' \subset J, \eta \in \hat{X}(T_{J'}), w, v \in W \hbox{
such that }
\\ &w \d(\xi)=w' \d(\eta), v \xi=v' \eta\} \hbox{ and } [J', w', v'] \\
& \hbox{ is contained in the closure of } [J, w, v] \Bigr\}.
\end{align*}

Then we have that $$<A_{\eta, J', w', v'}>_{(J', w', v', \eta) \in
\ci_{J, w, v, \xi}}=<(i_{J', w', v'})_! \cl_{\eta, J', w',
v'}>_{(J', w', v', \eta) \in \ci_{J, w, v, \xi}},$$ where $i_{J',
w', v'}: [J', w', v'] \hookrightarrow \overline{G^1}$ is the
inclusion.
\end{th3}

\subsection*{3.3} For $J \subset I$, set $B_J=(B \cap L_J)/Z^0(L_J)$. For $x \in W$,
set $G_{J, x}=B \dot x B/Z^0(L_J)$ and
$G'_{J, x}=B \dot x U_{P^-_{\d(J)}} (B \cap
L_{\d(J)})/Z^0(L_{\d(J)})$.

Define $\phi_{J, x}: G_{J, x} \rightarrow T_J$ by $\phi_{J, x} (u
\dot x t u')=t$ for $t \in T_J$ and $u, u' \in U$. Define
$\phi'_{J, x}: G'_{J, x} \rightarrow T_{\d(J)}$ by $\phi_{J, x} (u
\dot x t u')=t$ for $t \in T_{\d(J)}$ and $u \in U$. and $u' \in
U_{P^-_{\d(J)}} (U \cap L_{\d(J)})$.

Let $\xi \in \hat{X}(T_J)$. Then $\phi_{J, x}^* \cl_{\xi}=\cl_{J,
x, \xi}$ is a (tame) local system on $G_{J, x}$ and $(\phi'_{J,
x})^* \cl_{\d(\xi)}=\cl'_{J, x, \xi}$ is a (tame) local system on
$G'_{J, x}$. We denote by $A_{J, x, \xi}$ the perverse extension
of $\cl_{J, x, \xi}$ to $G/Z^0(L_J)$ and $A'_{J, x, \xi}$ the
perverse extension of $\cl'_{J, x, \xi}$ to $G/Z^0(L_{\d(J)})$. We
simply write $\cl_{I, x, \xi}$ as $\cl_{x, \xi}$ and $A_{I, x,
\xi}$ as $A_{x, \xi}$.

By \cite[4.1.2]{MS1}, $\cl_{J, x, \xi}$ and $A_{J, x, \xi}$ have
weight $x \xi$ for the left $B$-action and $-\xi$ for the right
$B_J$-action. Similarly, $\cl'_{J, x, \xi}$ and $A'_{J, x, \xi}$
have weight $x \d(\xi)$ for the left $B$-action and $-\d(\xi)$ for
the right $B_{\d(J)}$-action.

\subsection*{3.4} The group $B$ acts on $G \times G/Z^0(L_J)$ by $b
(g, h)=(g b \i, b h)$. A quotient $G \times^B G/Z^0(L_J)$ exists
and the product map $G \times G/Z^0(L_J) \rightarrow G/Z^0(L_J)$
defined by $(g, h) \mapsto g h$ induces a proper morphism $$m_J: G
\times^B G/Z^0(L_J) \rightarrow G/Z^0(L_J).$$

It is easy to see that $(m_J)_! (A_{w, \xi} \odot A_{J, 1,
\xi})=A_{J, w, \xi}$.

For $\xi \in \hat{X}(T_J)$ and $w, x \in W$, set $$\ci(\xi, J, w,
x)=\{w' \mid w' \xi=w x \xi \hbox{ and } G_{J, w'} \subset
\overline{m_J(G_{I, w} \times G_{J, x})}\},$$ where
$\overline{m_J(G_{I, w} \times G_{J, x})}$ is the closure of
$m_J(G_{I, w} \times G_{J, x})$ in $G/Z^0(L_J)$. Then by
\cite[4.2]{MS1}, $$(m_J)_! (A_{w, x \xi} \odot A_{J, x, \xi}) \in
<A_{J, w', \xi}>_{w' \in \ci(\xi, J, w, x)}.$$

Similarly,  set $$\ci'(\xi, J, w, x)=\{w' \mid w' \d(\xi)=w x
\d(\xi) \hbox{ and } G'_{J, w'} \subset \overline{m_J(G_{I, w}
\times G'_{J, x})}\}.$$ Then $$(m_J)_! (A_{w, x \xi} \odot A'_{J,
x, \xi}) \in <A'_{J, w', \xi}>_{w' \in \ci'(\xi, J, w, x)}.$$

\subsection*{3.5} The group $B_J$ acts on $G/Z^0(L_{\d(J)}) \times G/Z^0(L_J)$ by $b
(g, g')=(g g_0 b \i g_0 \i, g' b \i)$. The quotient
$G/Z^0(L_{\d(J)}) \times^{B_J} G/Z^0(L_J)$ exists and the map $G
\times G \rightarrow \overline{G^1}$ defined by $(g, g') \mapsto
(g, g') \cdot h_J$ induces a morphism $$p_J: G/Z^0(L_{\d(J)})
\times^{B_J} G/Z^0(L_J) \rightarrow \overline{G^1}.$$

For $w \in W^{\d(J)}$ and $v \in W$, we have that \begin{align*}
(B \dot w U_{P^-_{\d(J)}} (B \cap L_{\d(J)}), B \dot v B) \cdot
h_J &=(B \dot w (B \cap L_{\d(J)}), B \dot v (B \cap L_J))
\cdot h_J \\ &=(B \dot w (B \cap L_{\d(J)}), B \dot v) \cdot h_J \\
&=(B \dot w, B \dot v) \cdot h_J=[J, w, v]. \end{align*}

We can see that $p_J \mid_{G'_{J, w} \times^{B_J} G_{J, v}}$
factors through an affine space bundle map $p_{J, w, v}: G'_{J, w}
\times^{B_J} G_{J, v} \rightarrow [J, w, v]$. For $\xi \in
\hat{X}(T_J)$, $\cl'_{J, w, \xi} \odot \cl_{J, v, -\xi}$ is a
local system on $G'_{J, w} \times^{B_J} G_{J, v}$ and $$(p_J
\mid_{G'_{J, w} \times^{B_J} G_{J, v}})_! (\cl'_{J, w, \xi} \odot
\cl_{J, v, -\xi})=(i_{J, w, v})_! \cl_{\xi, J, w, v}.$$

For $w \in W^{\d(J)}$ and $u \in W_J$, $p_J(G'_{J, w \d(u)}
\times^{B_J} G_{J, 1})=[J, w, u \i]$ and $p_J \mid_{G'_{J, w
\d(u)} \times^{B_J} G_{J, 1}}$ factors through an affine space
bundle map $p_{J, w \d(u), 1}: G'_{J, w \d(u)} \times^{B_J} G_{J,
1} \rightarrow [J, w, u \i]$. For $\xi \in \hat{X}(T_J)$, we have
that \begin{align*} (i_{J, w, u \i})_! \cl_{u \xi, J, w, u \i}
&=(p_J \mid_{G'_{J, w \d(u)} \times^{B_J} G_{J, 1}})_! (\cl'_{J, w
\d(u), \xi} \odot \cl_{J, 1, -\xi}) \\ &=(p_J \mid_{G'_{J, w}
\times^{B_J} G_{J, u \i}})_! (\cl'_{J, w, u \xi} \odot \cl_{J, u
\i, -u \xi}).\end{align*}

\subsection*{3.6} The group $B$ acts on $G \times \overline{G^1}$ by $b(g, z)=(g b
\i, (b, b) z)$ for $b \in B, g \in G$ and $z \in \overline{G^1}$.
The quotient $G \times_B \overline{G^1}$ exists. The map $G \times
\overline{G^1} \rightarrow \overline{G^1}$ defined by $(g, z)
\mapsto (g, g) \cdot z$ induces a morphism $\rho: G \times_B
\overline{G^1} \rightarrow \overline{G^1}$. The map $\rho$ is the
unique extension of the map $\g$ defined in 2.3.

For $J \subset I$ and $w, v \in W$, we denote by $\rho_{J, w, v}$
the restriction of $\rho$ to $G \times_B [J, w, v]$.

Let $\xi \in \hat{X}(T_J) \subset \hat{X}$ with $w \d(\xi)=v \xi$.
Then $\bar{\mathbb Q}_l \odot \cl_{\xi, J, w, v}$ is a local
system on $G \times_B [J, w, v]$ and $\bar{\mathbb Q}_l[\dim(G)]
\odot A_{\xi, J, w, v}$ is a simple perverse sheaf on $G \times_B
\overline{G^1}$. Set $L_{\xi, J, w, v}=(\rho \mid_{G \times_B [J,
w, v]})_! (\bar{\mathbb Q}_l \odot \cl_{\xi, J, w, v})$ and
$C_{\xi, J, w, v}=\rho_! (\bar{\mathbb Q}_l \odot A_{\xi, J, w,
v})$. Then $L_{\xi, J, w, v}, C_{\xi, J, w, v} \in
\cd(\overline{G^1})$. Moreover, by the decomposition theorem in
\cite{BBD}, $C_{\xi, J, w, v}$ is semisimple.

The following result is an easy consequence of 3.2.

\begin{th4} We keep the notation of 3.2. Then $$<C_{\eta, J', w', v'}>_{(J', w', v', \eta) \in
\ci_{J, w, v, \xi}}=<L_{\eta, J', w', v'}>_{(J', w', v', \eta) \in
\ci_{J, w, v, \xi}}.$$
\end{th4}

\subsection*{3.8} Similarly, $B$ acts on $G \times (G/Z^0(L_{\d(J)}) \times^{B_J}
G/Z^0(L_J))$ by $b (g, z)=(g b \i, b \cdot z)$ for $b \in B$, $g
\in G$ and $z \in G/Z^0(L_{\d(J)}) \times^{B_J} G/Z^0(L_J)$, where
$B$ acts diagonally on $G/Z^0(L_{\d(J)}) \times^{B_J} G/Z^0(L_J)$
on both factors on the left. Denote by $G \times_B
(G/Z^0(L_{\d(J)}) \times^{B_J} G/Z^0(L_J))$ the quotient space.

The morphism $G \times G \times G \rightarrow \overline{G^1}$
defined by $(g, g_1, g_2) \mapsto (g g_1, g g_2) \cdot h_J$
induces a morphism $\pi_J: G \times_B (G/Z^0(L_{\d(J)})
\times^{B_J} G/Z^0(L_J)) \rightarrow \overline{G^1}$. We simply
write $\pi_J \mid_{G \times_B (G'_{J, w} \times^{B_J} G_{J, v})}$
as $\pi_{J, w, v}$. For $\xi \in \hat{X}(T_J)$ with $w \d(\xi)=v
\xi$, $\bar{\mathbb Q}_l \odot (\cl'_{J, w, \xi} \odot \cl_{J, v,
-\xi})$ is a local system on $G \times_B (G'_{J, w} \times^{B_J}
G_{J, v})$ and $$(\pi_{J, w, v})_! \bigl(\bar{\mathbb Q}_l \odot
(\cl'_{J, w, \xi} \odot \cl_{J, v, -\xi}) \bigr)=L_{\xi, J, w,
v}.$$

\begin{th5} Let $J \subset I$, $w, v \in W$ and $\xi \in
\hat{X}(T_J)$ with $w \d(\xi)=v \xi$. Then $$(\pi_J)_!
\bigl(\bar{\mathbb Q}_l \odot (A'_{J, w, \xi} \odot A_{J, v,
-\xi}) \bigr) \in <(\pi_J)_! \bigl(\bar{\mathbb Q}_l \odot (A'_{J,
w', \xi} \odot A_{J, 1, -\xi}) \bigr)>_{w' \in \ci'(\xi, J, v \i,
w)}.$$
\end{th5}

Proof. The group $B_J$ acts on $G \times G/Z^0(L_J)$ on the second
factor on the right. This action induces a $B_J$ action on $G
\times^B G/Z^0(L_J)$. Define the $B_J$ action on $G/Z^0(L_{\d(J)})
\times (G \times^B G/Z^0(L_J))$ by $b (g, z)=(g g_0 b \i g_0 \i, b
z)$ for $g \in G/Z^0(L_{\d(J)})$ and $z \in G \times^B
G/Z^0(L_J)$. The quotient $G/Z^0(L_{\d(J)}) \times^{B_J} (G
\times^B G/Z^0(L_J))$ exists and the morphism $(id, m_J):
G/Z^0(L_{\d(J)}) \times (G \times^B G/Z^0(L_J)) \rightarrow
G/Z^0(L_{\d(J)}) \times G/Z^0(L_J)$ induces a morphism
$$G/Z^0(L_{\d(J)}) \times^{B_J} (G \times^B G/Z^0(L_J)) \rightarrow
G/Z^0(L_{\d(J)}) \times^{B_J} G/Z^0(L_J).$$

The group $B$ action on $G \times G/Z^0(L_J)$ by acting on the
first factor on the left induces a $B$ action on $G \times^B
G/Z^0(L_J)$. Then we obtain the diagonal $B$-action on
$G/Z^0(L_{\d(J)}) \times^{B_J} (G \times^B G/Z^0(L_J))$. We may
also define the $B$-action on $G \times \bigl(G/Z^0(L_{\d(J)})
\times^{B_J} (G \times^B G/Z^0(L_J)) \bigr)$ in the same way as we
did in 3.8. We write the quotient space as $G \times_B
\bigl(G/Z^0(L_{\d(J)}) \times^{B_J} (G \times^B G/Z^0(L_J))
\bigr)$. The morphism $(id, id, m_J): G \times G/Z^0(L_{\d(J)})
\times (G \times^B G/Z^0(L_J)) \rightarrow G \times
G/Z^0(L_{\d(J)}) \times G/Z^0(L_J)$ induces a morphism
\begin{align*} m_{1, 2, 34}:{} & G \times_B \bigl(G/Z^0(L_{\d(J)})
\times^{B_J} (G \times^B G/Z^0(L_J)) \bigr) \rightarrow \\ & G
\times_B (G/Z^0(L_{\d(J)}) \times^{B_J} G/Z^0(L_J)).\end{align*}

The morphism $G \times G/Z^0(L_{\d(J)}) \times G \times G/Z^0(L_J)
\rightarrow \overline{G^1}$ defined by $(g_1, g_2, g_3, g_4)
\mapsto (g_1 g_2, g_1 g_3 g_4) \cdot h_J$ induces a morphism
$$f_{1, 2, 34}: G \times_B \bigl(G/Z^0(L_{\d(J)}) \times^{B_J} (G
\times^B G/Z^0(L_J)) \bigr) \rightarrow \overline{G^1}$$ and we
have that $f_{1, 2, 34}=\pi_J \circ m_{1, 2, 34}$. We also have
that $$(m_{1, 2, 34})_! \Bigl(\bar{\mathbb Q}_l \odot \bigl(A'_{J,
w, \xi} \odot (A_{v, -\xi} \odot A_{J, 1, -\xi}) \bigr)
\Bigr)=\bar{\mathbb Q}_l \odot (A'_{J, w, \xi} \odot A_{J, v,
-\xi})[\cdot].$$ Hence $$(\pi_J)_! \bigl(\bar{\mathbb Q}_l \odot
(A'_{J, w, \xi} \odot A_{J, v, -\xi}) \bigr)=(f_{1, 2, 34})_!
\Bigl(\bar{\mathbb Q}_l \odot \bigl(A'_{J, w, \xi} \odot (A_{v,
-\xi} \odot A_{J, 1, -\xi}) \bigr) \Bigr)[\cdot].$$

Similarly, we may define $G \times_B \bigl( (G \times^B
G/Z^0(L_{\d(J)})) \times^{B_J} G/Z^0(L_J) \bigr)$ and the morphism
$(id, m_{\d(J)}, id): G \times (G \times^B G/Z^0(L_{\d(J)}))
\times G/Z^0(L_J) \rightarrow G \times G/Z^0(L_{\d(J)}) \times
G/Z^0(L_J)$ induces a morphism \begin{align*} m_{1, 23, 4}:{} & G
\times_B \bigl( (G \times^B G/Z^0(L_{\d(J)})) \times^{B_J}
G/Z^0(L_J) \bigr) \rightarrow \\ & G \times_B (G/Z^0(L_{\d(J)})
\times^{B_J} G/Z^0(L_J)).\end{align*}

The morphism $G \times G \times G/Z^0(L_{\d(J)}) \times G/Z^0(L_J)
\rightarrow \overline{G^1}$ defined by $(g_1, g_2, g_3, g_4)
\mapsto (g_1 g_2 g_3, g_1 g_4) \cdot h_J$ induces a morphism
$$f_{1, 23, 4}: G \times_B \bigl( (G \times^B G/Z^0(L_{\d(J)})) \times^{B_J} G/Z^0(L_J)
\bigr) \rightarrow \overline{G^1}$$ and we have that $f_{1, 23,
4}=\pi_J \circ m_{1, 23, 4}$.

The isomorphism $G \times G/Z^0(L_{\d(J)}) \times G \times
G/Z^0(L_J) \rightarrow G \times G \times G/Z^0(L_{\d(J)}) \times
G/Z^0(L_J)$ defined by $(g_1, g_2, g_3, g_4) \mapsto (g_1 g_3, g_3
\i, g_2, g_4)$ induces an isomorphism \begin{align*}\iota:{} & G
\times_B \bigl(G/Z^0(L_{\d(J)}) \times^{B_J} (G \times^B
G/Z^0(L_J)) \bigr) \xrightarrow \simeq \\ & G \times_B \bigl( (G
\times^B G/Z^0(L_{\d(J)})) \times^{B_J} G/Z^0(L_J)
\bigr).\end{align*}

It is easy to see that $f_{1, 2, 34}=f_{1, 23, 4} \circ \iota$ and
$$\iota_! \Bigl(\bar{\mathbb Q}_l \odot \bigl(A'_{J, w, \xi} \odot
(A_{v, -\xi} \odot A_{J, 1, -\xi}) \bigr)
\Bigr)=\Bigl(\bar{\mathbb Q}_l \odot \bigl( (A_{v \i, v \xi} \odot
A'_{J, w, \xi}) \odot A_{J, 1, -\xi} \bigr) \Bigr)[\cdot].$$
Therefore
\begin{align*} (\pi_J)_! \bigl(\bar{\mathbb Q}_l \odot (A'_{J, w, \xi}
\odot & A_{J, v, -\xi}) \bigr)=(f_{1, 23, 4})_! \Bigl(\bar{\mathbb
Q}_l \odot \bigl( (A_{v \xi, v \i} \odot A'_{J, w, \xi}) \odot
A_{J, 1, -\xi} \bigr) \Bigr)[\cdot] \\ & \quad =(\pi_J)_! (m_{1,
23, 4})_! \Bigl(\bar{\mathbb Q}_l \odot \bigl( (A_{v \xi, v \i}
\odot A'_{J, w, \xi}) \odot A_{J, 1, -\xi} \bigr) \Bigr)[\cdot] \\
& \quad \in <(\pi_J)_! \bigl(\bar{\mathbb Q}_l \odot (A'_{J, w',
\xi} \odot A_{J, 1, -\xi}) \bigr)>_{w' \in \ci'(\xi, J, v \i, w)}.
\end{align*}

\begin{th6} Let $J \subset I$, $w, v \in W$ and $\xi \in
\hat{X}(T_J)$ with $w \d(\xi)=v \xi$. Then $$L_{\xi, J, w, v} \in
<L_{\xi, J, w', 1}>_{w' \in \ci'(\xi, J, v \i, w)}.$$
\end{th6}

Proof. We have that $$(\pi_{J, w, v})_! \bigl(\bar{\mathbb Q}_l
\odot (\cl'_{J, w, \xi} \odot \cl_{J, v, -\xi}) \bigr) \in
<(\pi_J)_! \bigl(\bar{\mathbb Q}_l \odot (A'_{J, w_1, \xi} \odot
A_{J, v_1, -\xi}) \bigr)>_{w_1, v_1},$$ where $w_1$ runs over the
elements in $W$ such that $G'_{J, w_1}$ is contained in the
closure of $G'_{J, w}$ and $v_1$ runs over the elements in $W$
such that $G_{J, v_1}$ is contained in the closure of $G_{J, v}$.

By lemma 3.9, we have that $$(\pi_J)_! \bigl(\bar{\mathbb Q}_l
\odot (A'_{J, w_1, \xi} \odot A_{J, v_1, -\xi}) \bigr) \in
<(\pi_J)_! \bigl(\bar{\mathbb Q}_l \odot (A'_{J, w_2, \xi} \odot
A_{J, 1, -\xi}) \bigr)>_{w_2},$$ where $w_2$ runs over the
elements in $W$ such that $w_2 \d(\xi)=\xi$ and $G'_{J, w_2}$ is
contained in the closure of $m_{\d(J)}(G_{I, v_1 \i} \times G'_{J,
w_1})$. Therefore,
\begin{align*} (\pi_{J, w, v})_! \bigl(\bar{\mathbb Q}_l \odot (\cl'_{J,
w, \xi} \odot \cl_{J, v, -\xi}) \bigr) & \in <(\pi_J)_!
\bigl(\bar{\mathbb Q}_l \odot
(A'_{J, w', \xi} \odot A_{J, 1, -\xi}) \bigr)>_{w'} \\
&=<(\pi_{J, w', 1})_! \bigl(\bar{\mathbb Q}_l \odot (\cl'_{J, w',
\xi} \odot \cl_{J, 1, -\xi}) \bigr)>_{w'},
\end{align*} where
$w'$ runs over the elements in $W$ such that $w' \d(\xi)=\xi$ and
$G'_{J, w'}$ is contained in the closure of $m_{\d(J)}(G_{I, v \i}
\times G'_{J, w})$. \qed

\begin{th7} Let $J \subset I$ and $\co$ be a $W_J$-orbit on
$\hat{X}(T_J)$. Set \begin{align*} \ci_{J, \co} &=\{(\xi, w, v)
\mid \xi \in \co, w, v \in W \hbox{ with } w \d(\xi)=v \xi\},
\\ \ci'_{J, \co} &=\{(\xi, v w, 1) \in \ci_{J, \co} \mid w \in
W^{\d(J)}, v \in W_{I(J, w, \d)}\}. \end{align*} Then
$$<L_{\xi, J, w, v}>_{(\xi, w, v) \in \ci_{J, \co}}=<L_{\xi, J, v w, 1
}>_{(\xi, v w, 1) \in \ci'_{J, \co}}.$$
\end{th7}

Proof. Set $J'=\{i \in I \mid \a_i=-w_0 \a_{\d(j)} \hbox{ for some
} j \in J\}$ and $y=w_0 w^{\d(J)}_0$. Then $y \d(J)=J'$.

Let $w \in W^{\d(J)}$ and $u \in W_{\d(J)}$. Then we have that $w
y \i \in W^{J'}$ and $l(w u y \i)=l(w y \i y u y \i)=l(w y \i)+l(y
u y \i)=l(w y \i)+l(u)$. Thus \begin{align*} \dim(G'_{J, w u})
&=\dim(G_{J, w u y \i})=l(w u y \i)+\dim(B/Z^0(L_J)) \\ &=l(w y
\i)+l(u)+\dim(B/Z^0(L_J)). \end{align*}

For $a, b \in W$, we have that $\dim(m_{\d(J)} (G_{I, a} \times
G'_{J, b})) \le l(a)+l(b y \i)+\dim(B/Z^0(L_J))$ and if the
equality holds, then $m_{\d(J)} (G_{I, a} \times G'_{J, b})=G'_{J,
a b}$.

Now by 3.10, $<L_{\xi, J, w, v}>_{(\xi, w, v) \in \ci_{J,
\co}}=<L_{\xi, J, w, 1}>_{(\xi, w, 1) \in \ci_{J, \co}}$.
Moreover, for $\xi \in \co$, $w \in W^{\d(J)}$ and $u \in W_J$
with $w \d(u) \d(\xi)=\xi$, we have that $L_{\xi, J, w \d(u),
1}=L_{u \xi, J, w, u \i} \in <L_{u \xi, J, w', 1}>_{w' \in \ci'(u
\xi, J, u, w)}$. By induction, it suffices to prove the following
statement:

Let $(w_i, u_i, \xi_i)_{i \ge 1}$ be a sequence with $w_i \in
W^{\d(J)}$, $u_i \in W_J$, $\xi_i \in \hat{X}(T_J)$. If $w_1
\d(u_1) \d(\xi_1)=\xi_1$ and for each $i \ge 1$, $\xi_{i+1}=u_i
\xi_i$ and $w_{i+1} \d(u_{i+1}) \in \ci'(\xi_{i+1}, J, u_i, w_i)$,
then for $n \gg 0$, $u_n \in W_{I(J, w_n, \d)}$.

For each $i \ge 1$, we have that
\begin{align*} l(w_{i+1} y \i)+l(u_{i+1}) &=\dim(G'_{J, w_{i+1}
\d(u_{i+1})})-\dim(B/Z^0(L_J))
\\ & \le \dim(m_{\d(J)} (G_{I, u_i} \times G'_{J,
w_i}))-\dim(B/Z^0(L_J))
\\ & \le l(w_i y \i)+l(u_i). \end{align*}

Moreover, if the equalities hold, then $G'_{J, w_{i+1}
\d(u_{i+1})}=G'_{J, u_i w_i}$ and $w_{i+1} \d(u_{i+1})=u_i w_i$.

Thus for $k \gg 0$, we have $l(w_k y \i)+l(u_k)=l(w_{k+1} y
\i)+l(u_{k+1})=\cdots$ and $w_{i+1} \d(u_{i+1})=u_i w_i$ for all
$i \ge k$. By \cite[3.10]{H2}, $w_{i+1} \ge w_i$ for $i \ge k$.
Therefore, for $m \gg k$, we have that
$w_m=w_{m+1}=\cdots=w_{\infty}$ and $l(u_m)=l(u_{m+1})=\cdots$.
Thus, $w_{\infty} \i u_i w_{\infty}=\d(u_{i+1}) \in W_{\d(J)}$ and
$l(w_{\infty} \i u_i w_{\infty})=l(u_{i+1})=l(u_i)$ for $i \ge m$.
So $w_{\infty} \i \supp(u_i) =\d(\supp(u_{i+1}))$ for all $i \ge
m$. In other words, for $i \ge m$, $$\cup_{l \ge i+1} \supp(u_l)
\subset \cup_{l \ge i} \supp(u_l)=w_{\infty} \d(\cup_{l \ge i+1}
\supp(u_l)).$$

Thus for $n \gg m$, $\cup_{i \ge n} \supp(u_i)=w_{\infty}
\d(\cup_{i \ge n} \supp(u_i))$. So $\supp(u_n) \subset \cup_{i \ge
l} \supp(u_i) \subset I(J, w_{\infty}, \d)$ and $u_n \in W_{I(J,
w_n, \d)}$. The lemma is proved. \qed

\

Combining 3.7 with 3.11, we obtain the key lemma.

\begin{th8} Let $\co$ be a $W$-orbit on $\hat{X}$. Set \begin{align*} \ci_{\co} &=\{(\xi, J, w, v)
\mid J \subset I, \xi \in \co \cap \hat{X}(T_J), w, v \in W \hbox{
with } w \d(\xi)=v \xi\}, \\ \ci'_{\co} &=\{(\xi, J, v w, 1) \in
\ci_{\co} \mid w \in W^{\d(J)}, v \in W_{I(J, w, \d)}\}.
\end{align*} Then
$$<C_{\xi, J, w, v}>_{(\xi, J, w, v) \in \ci_{\co}}=<L_{\xi, J,
w, v}>_{(\xi, J, w, v) \in \ci_{\co}}=<L_{\xi, J, v w, 1}>_{(\xi,
J, v w, 1) \in \ci'_{\co}}.$$
\end{th8}

\section{The character sheaves on $\overline{G^1}$}

By 2.1 and 3.12, we have the following result.

\begin{th9} Let $\co$ be a $W$-orbit on $\hat{X}$ and $A$ be a simple perverse sheaf
on $\overline{G^1}$. The following conditions on $A$ are
equivalent:

(i) $A \dashv C_{\xi, J, w, v}$ for some $J \subset I$, $w, v \in
W$ and $\xi \in \co \cap \hat{X}(T_J)$ with $w \d(\xi)=v \xi$.

(ii) $A \dashv L_{\xi, J, w, v}$ for some $J \subset I$, $w, v \in
W$ and $\xi \in \co \cap \hat{X}(T_J)$ with $w \d(\xi)=v \xi$.

(iii) $A \dashv L_{\xi, J, v w, 1}$ for some $J \subset I$, $w \in
W^{\d(J)}$, $v \in W_{I(J, w, \d)}$ and $\xi \in \co \cap
\hat{X}(T_J)$ with $v w \d(\xi)=\xi$.
\end{th9}

\subsection*{4.2} Let $\cc_{\co}(\overline{G^1})$ be the set of
(isomorphism classes) of simple perverse sheaves on
$\overline{G^1}$ which satisfies the equivalent conditions
4.1(i)-(iii) with respect to $\co$. The simple perverse sheaves on
$\overline{G^1}$ which belong to $\cc_{\co}(\overline{G^1})$ for
some $W$-orbit $\co$ are called character sheaves on the group
compactification $\overline{G^1}$; they (or their isomorphism
classes) form a set $\cc(\overline{G^1})$.

\subsection*{4.3} We keep the notation of 1.4. Now $K=I(J, w,
\d)$. Set $L=L_K/Z^0(L_J)$ and $B_1=(B \cap L_K)/Z^0(L_J)$. For $v
\in W_K$, set $L_v=B_1 \dot v \dot w g_0 B_1$.

Consider the following diagram
$$L \dot w g_0 \xleftarrow {p_1} G \times L \dot w g_0
\xrightarrow {p_2} G \times_{L_K} L \dot w g_0 \xrightarrow
{\pi^w_J} Z^w_J$$ where $p_1$ is the projection to the second
factor, $p_2$ is the projection map and $\pi^w_J$ is the map in
1.6.

For any character sheaf $X$ on $L_K \dot w g_0 /Z^0(L_J)$, let
$\tx$ be the simple perverse sheaf on $Z^w_J$ such that
$\tx=(\pi^w_J)_! (\bar{\mathbb Q}_l \odot X)[\cdot]$. (Then
$(\pi^w_J)^* \tx$ is a shift of $\bar{\mathbb Q}_l \odot X$. By
the commuting diagram in 1.5, the simple perverse sheaf $\tx$ on
$Z^w_J$ is the same as the simple perverse sheaf obtained in
\cite[11.12]{L3}.)

For $\xi \in \hat{X}(T_J)$, let $\cc_{\xi}(Z^w_J)$ be the
(isomorphism classes) of simple perverse sheaves on $Z^w_J$
consisting of all $\tx$ as above for $X \in \cc_{\xi}(L \dot w
g_0)$. The simple perverse sheaves on $Z^w_J$ which belong to
$\cc_{\xi}(Z^w_J)$ for some $\xi \in \hat{X}(T_J)$ are called
character sheaves on the $G$-stable piece $Z^w_J$; they (or their
isomorphism classes) form a set $\cc(Z^w_J)$.

\

Now we study $\cc_{\xi}(Z^w_J)$. For $v \in W_K$ with $v w
\d(\xi)=\xi$, let $\cl^L_{\xi, v, w \d}$ be the tame local system
on $L_v$ with weight $v w \d(\xi)$ for the left $B_1$-action and
$A^L_{\xi, v, w \d}$ be the perverse extension of $\cl^L_{\xi, v,
w \d}$ to $L$.

We identify $L \times_{B_1} L_v$ with $L_K \times_{B \cap L_K}
L_v$. Let $\bar{\mathbb Q}^G_l$ be the trivial local system on $G$
and $\bar{\mathbb Q}^{L_K}_l$ be the trivial local system on
$L_K$. Define $\rho^L_v: L_K \times_{B \cap L_K} L_v \rightarrow L
\dot w g_0$ by $\rho(l, l')=l l' l \i$. Then for $X \in
\cc_{\xi}(L \dot w g_0 )$, there exists $v \in W_K$ with $v w
\d(\xi)=\xi$ and $X \dashv (\rho^L_v)_! (\bar{\mathbb Q}^{L_K}_l
\odot A^L_{\xi, v, w \d})$. By the decomposition theorem of
\cite{BBD}, $(\rho^L_v)_! (\bar{\mathbb Q}^{L_K}_l \odot A^L_{\xi,
v, w \d})$ is semisimple. We have that $(\rho^L_v)_! (\bar{\mathbb
Q}^{L_K}_l \odot A^L_{\xi, v, w \d})=X[m] \oplus B$ for some $B
\in \cd(L)$ and $m \in \mathbb Z$. Therefore, $(\pi^w_J)_! \Bigl(
\bar{\mathbb Q}^G_l \odot \bigl( (\rho^L_v)_! (\bar{\mathbb
Q}^{L_K}_l \odot A^L_{\xi, v, w \d}) \bigr) \Bigr)=\tilde{X}[m']
\oplus C$ for some $C \in \cd(Z^w_J)$ and $m' \in \mathbb Z$. In
other words, a simple perverse sheaf $A$ on $Z^w_J$ is contained
in $\cc_{\xi}(Z^w_J)$ if and only if $A$ is an irreducible
constitute of $(\pi^w_J)_! \Bigl( \bar{\mathbb Q}^G_l \odot \bigl(
(\rho^L_v)_! (\bar{\mathbb Q}^{L_K}_l \odot A^L_{\xi, v, w \d})
\bigr) \Bigr)$ for some $v \in W_K$ with $v w \d(\xi)=\xi$.

Note that $<\cl^L_{\xi, v, w \d}>_{v \in W_K, v w
\d(\xi)=\xi}=<A^L_{\xi, v, w \d}>_{v \in W_K, v w \d(\xi)=\xi}$.
Then a simple perverse sheaf $A$ on $Z^w_J$ is contained in
$\cc_{\xi}(Z^w_J)$ if and only if $A \dashv (\pi^w_J)_! \Bigl(
\bar{\mathbb Q}^G_l \odot \bigl( (\rho^L_v)_! (\bar{\mathbb
Q}^{L_K}_l \odot \cl^L_{\xi, v, w \d}) \bigr) \Bigr)$ for some $v
\in W_K$ with $v w \d(\xi)=\xi$.

\subsection*{4.4} We keep the notation of 4.3. Consider the
commuting diagram

$$\begin{matrix} L_K \times_{B \cap L_K} L_v & \mapleft{p^L_1} & G \times
(L_K \times_{B \cap L_K} L_v) & \mapright{p^L_2} & G \times_{B \cap L_K} L_v \\
\mapdown{\rho^L_v} & & \mapdown{(id, \rho^L_v)} & &
\mapdown{\pi^L_v} \\ L \dot w g_0 & \mapleft{p_1} & G \times L
\dot w g_0 & \mapright{p_2} & G \times_{L_K} L \dot w g_0,
\end{matrix}$$ where $p^L_1$ and $p^L_2$ are projection maps and $$\pi^L_v: G \times_{B_1} L_v=G
\times_L (L \times_{B_1} L_v) \rightarrow G \times_L L \dot w
g_0$$ is the map induced from $(id, \rho^L_v)$. Then
\begin{align*} p_2^* & \Bigl(\bar{\mathbb Q}^G_l \odot \bigl(
(\rho^L_v)_! (\bar{\mathbb Q}^{L_K}_l \odot \cl^L_{\xi, v, w \d})
\bigr) \Bigr) =p_1^* (\rho^L_v)_! (\bar{\mathbb Q}^{L_K}_l \odot
\cl^L_{\xi, v, w \d}) \\ &=(id, \rho^L_v)_! (p^L_1)^*
(\bar{\mathbb Q}^{L_K}_l \odot \cl^L_{\xi, v, w \d})=(id,
\rho^L_v)_! (p^L_2)^* (\bar{\mathbb Q}^G_l \odot \cl^L_{\xi, v, w
\d}) \\ &=p_2^* (\pi^L_v)_! (\bar{\mathbb Q}^G_l \odot \cl^L_{\xi,
v, w \d}).
\end{align*}

Therefore $\bar{\mathbb Q}^G_l \odot \bigl( (\rho^L_v)_!
(\bar{\mathbb Q}^{L_K}_l \odot \cl^L_{\xi, v, w \d})
\bigr)=(\pi^L_v)_! (\bar{\mathbb Q}^G_l \odot \cl^L_{\xi, v, w
\d})$.

Define $\pi_{B_1}: B/Z^0(L_J) \rightarrow B_1$ by $\pi_{B_1}(u
b)=b$ for $u \in U_{P_{I(J, w, \d)}}$ and $b \in B_1$. Define the
map $B \dot v \dot w g_0 (B \cap L_J)/Z^0(L_J) \rightarrow L_v$ by
$b_1 \dot v \dot w g_0 b_2 \mapsto \pi_{B_1}(b_1) \dot v \dot w
g_0 \pi_{B_1}(b_2)$ for $b_1 \in B/Z^0(L_J)$ and $b_2 \in (B \cap
L_J)/Z^0(L_J)$. It is easy to see that this map is defined and is
an affine space bundle morphism. This map induces an affine space
morphism $$p^L_v: G \times_{B \cap L_K} B \dot x \dot w g_0 (B
\cap L_J)/Z^0(L_J) \rightarrow G \times_{B \cap L_K} L_v.$$

We also have $(p^L_v)_! (p^L_v)^* (\bar{\mathbb Q}^G_l \odot
\cl^L_{\xi, v, w \d})=(\bar{\mathbb Q}^G_l \odot \cl^L_{\xi, v, w
\d})[\cdot]$.

Set $\phi=\pi^w_J\circ \pi^L_v \circ p^L_v: G \times_{B_1} B \dot
x \dot w g_0 (B \cap L_J)/Z^0(L_J) \rightarrow Z^w_J$. Then
\begin{align*} (\pi^w_J)_! \Bigl( \bar{\mathbb Q}^G_l \odot \bigl(
(\rho^L_v)_! (\bar{\mathbb Q}^{L_K}_l \odot \cl^L_{\xi, v, w \d})
\bigr) \Bigr) &=(\pi^w_J)_! (\pi^L_v)_! (\bar{\mathbb Q}^G_l \odot
\cl^L_{\xi, v, w \d}) \\ &=\phi_! (p^L_v)^* (\bar{\mathbb Q}^G_l
\odot \cl^L_{\xi, v, w \d})[\cdot]. \end{align*}

Note that $\phi$ sends $(g, g')$ to $(g g' g_0 \i, g) \cdot h_J$
for $g \in G$ and $g' \in B \dot v \dot w g_0 (B \cap
L_J)/Z^0(L_J)$. Now we identify $B \dot v \dot w g_0 (B \cap
L_J)/Z^0(L_J)$ with $B_1 \times_{B_1} B \dot v \dot w g_0 (B \cap
L_J)/Z^0(L_J)$. Then \begin{align*} \phi(B \dot v \dot w g_0 (B
\cap L_J)/Z^0(L_J)) &=(B \dot v \dot w (B \cap L_{\d(J)}), 1)
\cdot h_J \\ &=(B \dot v \dot w, B \cap L_J) \cdot h_J=[J, v w,
1]. \end{align*}

By \cite[1.12]{H2}, $[J, v w, 1] \subset (P_K \dot w, P_K) \cdot
h_J \subset Z^w_J$. Therefore $$\rho(G \times_B [J, v w, 1])
\subset Z^w_J.$$ We denote by $\rho^{v w}_J: G \times_B [J, v w,
1] \rightarrow Z^w_J$ the restriction map of $\rho$.

Now the map $\phi: G \times_{B \cap L_K} G_{x w, \d} \rightarrow
Z^w_J$ factors through $$G \times_{B \cap L_K} B \dot v \dot w g_0
(B \cap L_J)/Z^0(L_J) \rightarrow G \times_{B \cap L_K} [J, v w,
1] \rightarrow G \times_B [J, v w, 1] \xrightarrow {\rho^{v w}_J}
Z^w_J,$$ where the first two maps are affine space bundle maps.
Hence
\begin{align*} (\pi^w_J)_! \Bigl( \bar{\mathbb Q}^G_l \odot \bigl(
(\rho^L_v)_! & (\bar{\mathbb Q}^{L_K}_l \odot \cl^L_{\xi, v, w
\d}) \bigr) \Bigr) =\phi_! (p^L_v)^* (\bar{\mathbb Q}^G_l \odot
\cl^L_{\xi, v, w \d})[\cdot] \\ &=(\rho^{v w}_J)_! (\bar{\mathbb
Q}^G_l \odot \cl_{\xi, J, v w, 1})[\cdot]=L_{\xi, J, v w, 1}
\mid_{Z^w_J}[\cdot].
\end{align*}

Since $\rho(G \times_B [J, v w, 1]) \subset Z^w_J$, we have
$L_{\xi, J, v w, 1} \mid_{Z^{w'}_{J'}}=0$ if $J \neq J'$ or $w
\neq w'$.

In summary, for $(\xi, J, w, v) \in \ci'_{\co}$, $$L_{\xi, J, v w,
1}\mid_{Z^{w'}_{J'}}=\begin{cases} (\pi^w_J)_! \Bigl( \bar{\mathbb
Q}^G_l \odot \bigl( (\rho^L_v)_! (\bar{\mathbb Q}^{L_K}_l \odot
\cl^L_{\xi, v, w \d}) \bigr) \Bigr)[\cdot], & \hbox{ if } J=J'
\hbox{ and } w=w'; \\ 0, & \hbox{ otherwise }.
\end{cases}$$

\

Now we will show that the character sheaves on $\overline{G^1}$
have the following ``nice'' property.

\begin{th10} Let $J \subset I$, $w \in W^{\d(J)}$ and $\co$ be a $W$-orbit on $\hat{X}$.
Let $A \in \cc_{\co}(\overline{G^1})$ and $K$ be a simple perverse
sheaf on $Z^w_J$. If $K \dashv A \mid_{Z^w_J}$, then $K \in
\cc_{\xi}(Z^w_J)$ for some $\xi \in \co \cap \hat{X}(T_J)$.
\end{th10}

Proof. We keep the notation of 3.12.

There exists $(\xi, J', w', v) \in \ci_\co$, such that $A \dashv
C_{\xi, J', w', v}$. Since $C_{\xi, J', w', v}$ is semisimple, we
have that $C_{\xi, J', w', v}=A[m] \oplus B$ for some $B \in
\cd(\overline{G^1})$ and $m \in \mathbb Z$. Then $C_{\xi, J', w',
v} \mid_{Z^w_J}=A[m] \mid_{Z^w_J} \oplus B\mid_{Z^w_J}$. Thus $K
\dashv C_{\xi, J', w', v}\mid_{Z^w_J}$.

By 2.1 and 3.12, $<C_{\xi, J', w', v}\mid_{Z^w_J}>_{(\xi, J', w',
v) \in \ci_{\co}}=<L_{\xi, J', v w', 1}\mid_{Z^w_J}>_{(\xi, J',
w', v) \in \ci'_{\co}}$. So $K \dashv (\pi^w_J)_! \Bigl(
\bar{\mathbb Q}^G_l \odot \bigl( (\rho^L_v)_! (\bar{\mathbb
Q}^{L_{I(J, w, \d)}}_l \odot \cl^L_{\xi, v, w \d}) \bigr) \Bigr)$
for some $v \in W_{I(J, w, \d)}$ and $\xi \in \co \cap
\hat{X}(T_J)$ with $v w \d(\xi)=\xi$. By 4.3, $K \in
\cc_{\xi}(Z^w_J)$. The lemma is proved. \qed

\begin{th11} Let $A$ be a simple perverse sheaf on $\overline{G^1}$. Then
$A \in \cc(\overline{G^1})$ if and only if there exists $J \subset
I$, $w \in W^{\d(J)}$ and $X \in \cc(Z^w_J)$, such that $A$ is the
perverse extension of $X$.
\end{th11}

\begin{rmk} Therefore, our definition of character sheaves on
$\overline{G^1}$ coincides with Lusztig's definition in
\cite[12.3]{L3}.
\end{rmk}

Proof. We follow the proof of \cite[11.15 \& 11.18]{L3}.

Since $\overline{G^1}=\bigsqcup_{J \subset I} \bigsqcup_{w \in
W^{\d(J)}} Z^w_J$, there exists $J \subset I$ and $w \in
W^{\d(J)}$, such that $\supp(A) \cap Z^w_J$ is dense in
$\supp(A)$. Then $A \mid_{Z^w_J}$ is a simple perverse sheaf on
$Z^w_J$ and $A$ is the perverse extension of $A \mid_{Z^w_J}$. By
4.5, $A \mid_{Z^w_J} \in \cc(Z^w_J)$.

On the other hand, if $X \in \cc(Z^w_J)$, then by 4.4, we may
assume that $X \dashv (\rho^{v w}_J)_! (\bar{\mathbb Q}_l \odot
\cl_{\xi, J, v w, 1})$ for some $\xi \in \hat{X}(T_J)$ with $v w
\d(\xi)=\xi$. Then there exists a simple perverse sheaf $A$ on
$\overline{G^1}$ such that $\supp(A)$ is the closure of $\supp(X)$
in $\overline{G^1}$, $A \mid_{Z^w_J}=X$ and $A \dashv \rho_!
(\bar{\mathbb Q}_l \odot \cl_{\xi, J, v w, 1})$. Then $A$ is the
perverse extension of $X$ and $A \in \cc(\overline{G^1})$. \qed

\begin{th12} Let $\co_1, \co_2$ be two distinct $W$-orbits on $\hat{X}$.
Then for $(\xi_1, J_1, w_1, v_1) \in \ci_{\co_1}$ and $(\xi_2,
J_2, w_2, v_2) \in \ci_{\co_2}$, we have that
$$H^{\cdot}_c(\overline{G^1}, C_{\xi_1, J_1, w_1, v_1} \otimes
C_{-\xi_2, J_2, w_2, v_2})=0.$$
\end{th12}

Proof. We keep the notation of 4.3. For any $W_K$-orbit $\co$ on
$\hat{X}(T_J)$, set
$$\ci^L_{\co}=\{(v, \xi) \mid v \in W_K, \xi \in \co \hbox{ with }
v w \d(\xi)=\xi\}.$$

We have that $<\cl^L_{\xi, v, w \d}>_{(v, \xi) \in
\ci^L_{\co}}=<A^L_{\xi, v, w \d}>_{(v, \xi) \in \ci^L_{\co}}$. Now
set $C^L_{\xi, v, w \d}=(\pi^w_J)_! \Bigl( \bar{\mathbb Q}^G_l
\odot \bigl( (\rho^L_v)_! (\bar{\mathbb Q}^{L_K}_l \odot A^L_{\xi,
v, w \d}) \bigr) \Bigr)$. Then $C^L_{\xi, v, w \d}$ is semisimple
and all the irreducible constitutes are contained in
$\cc_{\xi}(Z^w_J)$. Moreover, $$<L_{\xi, J, v w, 1}
\mid_{Z^w_J}>_{(v, \xi) \in \ci^L_{\co}}=<C^L_{\xi, v, w \d}>_{(v,
\xi) \in \ci^L_{\co}}.$$

Let $\co, \co'$ be two distinct $W_K$-orbits on $\hat{X}(T_J)$.
Then by 2.6 and \cite[1.2.5]{MS1}, for $(v, \xi) \in \ci^L_{\co}$
and $(v', \xi') \in \ci^L_{\co'}$, we have that
$$H^{\cdot}_c(Z^w_J, C^L_{\xi, v, w \d} \otimes C^L_{-\xi', v', w
\d})=0.$$

By 2.1, we also have that $H^{\cdot}_c(Z^w_J, L_{\xi, J, v w, 1}
\mid_{Z^w_J} \otimes L_{-\xi', J, v' w, 1} \mid_{Z^w_J})=0$.

By 4.4, for $(\xi_1, J_1, v_1 w_1, 1) \in \ci'_{\co_1}$ and
$(\xi_2, J_2, v_2 w_2, 1) \in \ci'_{\co_2}$,
$$H^{\cdot}_c(Z^w_J, L_{\xi_1, J_1, v_1 w_1, 1} \mid_{Z^w_J} \otimes
L_{-\xi_2, J_2, v_2 w_2, 1} \mid_{Z^w_J})=0.$$

Since the above equality holds for all $G$-stable pieces $Z^w_J$
of $\overline{G^1}$, we have that $H^{\cdot}_c(\overline{G^1},
L_{\xi_1, J_1, v_1 w_1, 1} \otimes L_{-\xi_2, J_2, v_2 w_2,
1})=0$. Then by 2.1 and 3.12, $$H^{\cdot}_c(\overline{G^1},
C_{\xi_1, J_1, w_1, v_1} \otimes C_{-\xi_2, J_2, w_2, v_2})=0.$$

\subsection*{4.8} As in \cite[11.12(c)]{L1}, we have that
$\cc_{\co}(\overline{G^1}) \cap
\cc_{\co'}(\overline{G^1})=\varnothing$ for distinct $W$-orbits
$\co$ and $\co'$ on $\hat{X}$. In other words, there is a well
defined map $\cc(\overline{G^1}) \rightarrow {\hbox{$W$-orbit on }
\hat{X}}$ given by attaching $A \in \cc(\overline{G^1})$ the
$W$-orbit $\co$, where $A \in \cc_{\co}(\overline{G^1})$.

\section{The parabolic character sheaves}

In this section, $G$ is a connected reductive algebraic group. We
keep the notation in 1.3.

\subsection{} We first study the closures of the $B \times
B$-orbits in $\tz_{J, y, \d}$.

As in 3.5, we define the action of $B \cap L_J$ on $G/{^{\dot y
\i} U_{P_{J'}}} \times G/{U_{P_J}}$ by $b (g, g')=(g g_0 b \i g_0
\i, g' b \i)$. Denote by $G/{^{\dot y \i} U_{P_{J'}}} \times^{B
\cap L_J} G/{U_{P_J}}$ the quotient space. The morphism $G \times
G \rightarrow \tz_{J, y, \d}$ defined by $(g_1, g_2) \mapsto (g_1,
g_2) \cdot \tilde h_{J, y, \d}$ induces a proper morphism $$p_{J,
y, \d}: G/{^{\dot y \i} U_{P_{J'}}} \times^{B \cap L_J}
G/{U_{P_J}} \rightarrow \tz_{J, y, \d}.$$

By 3.5, for $w \in W^{\d(J)}$ and $v \in W$, we have that
$$p_{J, y, \d} \bigl( (B \dot w \dot y \i B \dot y)/{^{\dot y \i}
U_{P_{J'}}} \times^{B \cap L_J} (B \dot v B)/{U_{P_J}}
\bigr)=\widetilde{[J, w, v]}_{y, \d}.$$

Moreover, for $w \in W^{\d(J)}$, $u \in W_J$ and $v \in W$. we
have that \begin{align*} (B \dot w \dot{\d(u)} \dot y \i B \dot y,
B \dot v B) \cdot \tilde h_{J, y, \d} &=\bigl(B \dot w \dot{\d(u)}
(B \cap L_{\d(J)}), B \dot v B \bigr) \cdot \tilde h_{J, y, \d} \\
&=(B \dot w, B \dot v B \dot u \i) \cdot \tilde h_{J, y, \d}
\\ & \subset \bigcup_{v_1 \le v, u_1 \le u} (B \dot w, B \dot v_1 \dot u_1 \i B)
\cdot \tilde h_{J, y, \d} \\ &=\bigcup_{v_1 \le v, u_1 \le u}
\widetilde{[J, w, v_1 u_1 \i]}_{y, \d}.\end{align*}

\begin{th13} Let $w \in W^{\d(J)}$ and $v \in W$. Denote by
$\overline{[J, w, v]}_{y, \d}$ the closure of $\widetilde{[J, w,
v]}_{y, \d}$ in $\tz_{J, y, \d}$. Then $$\overline{[J, w, v]}_{y,
\d}=\bigsqcup_{w_1 \in W^{\d(J)}, u_1 \in W_J, w_1 \d(u_1) y \i
\le w y \i} \bigsqcup_{v' \le v} \widetilde{[J, w_1, v' u_1
\i]}_{y, \d}.$$
\end{th13}

Proof. Note that $\bigsqcup_{w' \in W, w' y \i \le w y \i} B \dot
w' \dot y \i B \dot y$ and $\bigsqcup_{v' \le v} B \dot v' B$ are
irreducible closed subvarieties of $G$ and $p_{J, y, \d}$ is a
proper map. Thus we have that $$\overline{[J, w, v]}_{y,
\d}=\bigcup_{w' y \i \le w y \i, v' \le v} p_{J, y, \d} \bigl( (B
\dot w' \dot y \i B \dot y)/{^{\dot y \i} U_{P_{J'}}} \times^{B
\cap L_J} (B \dot v' B)/{U_{P_J}} \bigr).$$

Let $w_1 \in W^{\d(J)}$ and $u_1 \in W_J$. Then we have that $w_1
y \i \in W^{J'}$ and $y \d(u_1) y \i \in W_{J'}$. Moreover, if
$u_2 \le u_1$, then $y \d(u_2) y \i \le y \d(u_1) y \i$. Thus $w_1
\d(u_2) y \i=w_1 y \i y \d(u_1) y \i \le w_1 y \i y \d(u_1) y
\i=w_1 \d(u_1) y \i$. Therefore for $w_1 \d(u_1) y \i \le w y \i$
and $v' \le v$, we have that \begin{align*} & p_{J, y, \d} \bigl(
(B \dot w_1 \dot{\d(u_1)} \dot y \i B \dot y)/{^{\dot y \i}
U_{P_{J'}}} \times^{B \cap L_J} (B \dot v' B)/{U_{P_J}} \bigr) \\
& \subset \bigcup_{u_2 \in W_J, w_1 \d(u_2) y \i \le w_1 \d(u_1) y
\i \le w y \i} \bigcup_{v_1 \le v' \le v} \widetilde{[J, w_1, v_1
u_2 \i]}_{y, \d}.\end{align*}

Hence $\overline{[J, w, v]}_{y, \d} \subset \bigsqcup_{w_1 \in
W^{\d(J)}, u_1 \in W_J, w_1 \d(u_1) y \i \le w y \i} \bigsqcup_{v'
\le v} \widetilde{[J, w_1, v' u_1 \i]}_{y, \d}$.

On the other hand, for $w_1 \in W^{\d(J)}, u_1 \in W_J$ with $w_1
\d(u_1) y \i \le w y \i$ and $v' \le v$, we have that
\begin{align*} \widetilde{[J, w_1, v' u_1 \i]}_{y, \d} & =(B \dot
w_1, B \dot v' \dot u_1 \i) \cdot \tilde h_{J, y, \d}=(B \dot w_1
\dot{\d(u_1)}, B \dot v') \cdot \tilde h_{J, y, \d} \\ & \subset
p_{J, y, \d} \bigl( (B \dot w_1 \dot{\d(u_1)} \dot y \i B \dot
y)/{^{\dot y \i} U_{P_{J'}}} \times^{B \cap L_J} (B \dot v'
B)/{U_{P_J}} \bigr) \\ & \subset \overline{[J, w, v]}_{y, \d}.
\end{align*}

The lemma is proved. \qed

\subsection*{5.3} As a consequence of the description of the closure of $B \times
B$-orbits, we will describe the closure of the $G$-stable pieces
in $\tz_{J, y, \d}$ in 5.5. Although the result is not used in the
proofs of properties of parabolic character sheaves we are going
to discuss about, it serves as the motivation for them.

We have a bijection $\d_y=y \d: J \rightarrow J'$. Let $w, w' \in
W$ with $l(w)=l(w')$. We say that $w'$ can be obtained from $w$
via a $(J, \d_y)$-cyclic shift if $w=s_{i_1} s_{i_2} \cdots
s_{i_n}$ is a reduced expression and either (1) $i_1 \in J$ and
$w'=s_{i_1} w s_{\d_y(i_1)}$ or (2) $i_n \in J'$ and $w'=s_{\d_y
\i(i_n)} w s_{i_n}$. We write $w \sim_{J, \d_y} w'$ if there
exists a finite sequence of elements $w=w_0, w_1, \ldots, w_m=w'$
such that $w_{k+1}$ can be obtained from $w_k$ via a $(J,
\d_y)$-cyclic shift.

By \cite[3.9]{H2}, we have the following result.

\begin{th14} Let $w, w' \in W^{J'}$. Then the following conditions
are equivalent:

(1) $w \ge u w' \d_y(u) \i$ for some $u \in W_J$.

(2) $w \ge v w' \d_y(u) \i$ for some $v \le u \in W_J$.

(3) $w \ge x$ for some $x \sim_{J, \d_y} w'$.

In this case, we say that $w \ge_{J, \d_y} w'$.
\end{th14}

\begin{th15} Let $w \in W^{\d(J)}$. Then the closure of $\tz^w_{J,
y, \d}$ in $\tz_{J, y, \d}$ is $\bigsqcup_{w' \in W^{\d(J)}, w y
\i \ge_{J, \d_y} w' y \i} \tz^{w'}_{J, y, \d}$.
\end{th15}

The case where $y=1$ was proved in \cite[4.6]{H2}. The general
case can be treated in a similar way.

\

From now on, we study the parabolic character sheaves on $\tz_{J,
y, \d}$.

\subsection*{5.6} Let $w, v \in W$. Let $p_{w, v}: B
\times B \rightarrow \widetilde{[J, w, v]}_{y, \d}$ be the
morphism defined by $(b_1, b_2) \mapsto (b_1 \dot w, b_2 \dot v)
\cdot \tilde{h}_{J, y, \d}$ for $b_1, b_2 \in B$. By 1.3, $p_{w u,
v \d \i(u)}=p_{w, v}$ for $u \in W_{\d(J)}$.

The morphism $B \times B \rightarrow T$, $(t_1 u_1, t_2 u_2)
\mapsto (\dot w g_0) \i t_1 (\dot w g_0) \dot v \i t_2 \i \dot v$
factors through a morphism $pr_{w, v}: \widetilde{[J, w, v]}_{y,
\d} \rightarrow T$.

For $\xi \in \hat{X}$, $pr_{w, v}^* \cl_{\xi}=\cl_{\xi, w, v}$ is
a (tame) local system on $\widetilde{[J, w, v]}_{y, \d}$. We
denote by $A_{\xi, w, v}$ the perverse extension of $\cl_{\xi, w,
v}$ to $\tz_{J, y, \d}$.

\subsection*{5.7} We keep the notation of 1.4. Let $N_{\hat{G}}(L_K)$
be the normalizer of $L_K$ in $\hat{G}$. Note that
$N_{\hat{G}}(L_K)$ is a disconnected group with identity component
$L_K$ and $L_K \dot w g_0$ is a connected component of
$N_{\hat{G}}(L_K)$. Consider the diagram

$$L_K \dot w g_0 \xleftarrow {p_1} G \times L_K \dot w g_0
\xrightarrow {p_2} G \times_{L_K} L_K \dot w g_0 \xrightarrow
{\pi^w_{J, y, \d}} \tz^w_{J, y, \d}$$ where $p_1$ is the
projection to the second factor and $p_2$ is the projection map.

For any character sheaf $X$ on $L_K \dot w g_0$, let $\tx$ be the
simple perverse sheaf on $\tz^w_{J, y, \d}$ such that $\tx=(\tilde
\pi^w_{J, y, \d})_! (\bar{\mathbb Q}_l \odot X)[\cdot]$. Let
$\cc(\tz^w_{J, y, \d})$ be the (isomorphism classes) of simple
perverse sheaves on $\tz^w_{J, y, \d}$ consisting of all $\tx$ as
above. The elements in $\cc(\tz^w_{J, y, \d})$ are called the
character sheaves on $\tz^w_{J, y, \d}$.

\subsection*{5.8} The group $B$ acts on $G \times \tz_{J, y, \d}$ by $b(g,
z)=(g b \i, (b, b) z)$. The quotient $G \times_B \tz_{J, y, \d}$
exists. The map $G \times \tz_{J, y, \d} \rightarrow \tz_{J, y,
\d}$ defined by $(g, z) \mapsto (g, g) \cdot z$ induces a morphism
$\rho: G \times_B \tz_{J, y, \d} \rightarrow \tz_{J, y, \d}$.

By what we did in section 3 and section 4, we obtain the following
result.

\begin{th16} Let $A$ be a simple perverse sheaf on $\tz_{J, y, \d}$. The
following conditions are equivalent:

(i) $A \dashv \rho_! (\bar{\mathbb Q}_l \odot A_{\xi, w, v})$ for
$w, v \in W$ and $\xi \in \hat{X}$ with $w \d(\xi)=v \xi$.

(ii) $A \dashv (\rho \mid_{G \times_B \widetilde{[J, w, v]}_{y,
\d}})_! (\bar{\mathbb Q}_l \odot \cl_{\xi, w, v})$ for some $w, v
\in W$ and $\xi \in \hat{X}$ with $w \d(\xi)=v \xi$.

(iii) $A \dashv (\rho \mid_{G \times_B \widetilde{[J, v w, 1]}_{y,
\d}})_! (\bar{\mathbb Q}_l \odot \cl_{\xi, v w, 1})$ for some $w
\in W^{\d(J)}$, $v \in W_{I(J, w, \d)}$ and $\xi \in \hat{X}$ with
$v w \d(\xi)=\xi$.

(iv) There exists $w \in W^{\d(J)}$ and $X \in \cc(\tz^w_{J, y,
\d})$ such that $A$ is the perverse extension of $X$ to $\tz_{J,
y, \d}$.

Moreover, if $A$ satisfies the equivalent conditions, then for any
$w \in W^{\d(J)}$, any composition factor of $^p H^i(A
\mid_{\tz^w_{J, y, \d}})$ is contained in $\cc(\tz^w_{J, y, \d})$.
\end{th16}

\begin{rmk} The condition (iv) was one of the two equivalent
definitions of parabolic character sheaves in \cite{L3}. The fact
the $(i)\Leftrightarrow(iv)$ was first proved by Springer in
\cite{S2} in terms of $(P, Q, \s)$-character sheaves on $G$. The
``moreover'' part was first proved by Lusztig in \cite[11.14]{L3}.
Our approach gave a new proof of their results.
\end{rmk}

\subsection*{5.10} Let $\cc(\tz_{J, y, \d})$ be the (isomorphism
classes) of  simple perverse sheaves on $\tz_{J, y, \d}$ which
satisfy the equivalent conditions 5.9(i)-(iv). The elements of
$\cc(\tz_{J, y, \d})$ are called the parabolic character sheaves
on $\tz_{J, y, \d}$.

For $W_J$-orbit $\co$ on $\hat{X}$, we denote by
$\cc_{\co}(\tz_{J, y, \d})$ the (isomorphism classes) of simple
perverse sheaves $A$ on $\tz_{J, y, \d}$ such that $A$ satisfies
the condition 5.9(iii) for some $\xi \in \co$. Then as in 4.8, we
have that $$\cc_{\co}(\tz_{J, y, \d}) \cap \cc_{\co'}(\tz_{J, y,
\d})=\varnothing$$ for distinct $W_J$-orbits $\co$ and $\co'$ on
$\hat{X}$. In other words, there is a well defined map
$\cc(\tz_{J, y, \d}) \rightarrow {\hbox{$W_J$-orbit on } \hat{X}}$
given by attaching $A \in \cc(\tz_{J, y, \d})$ the $W$-orbit
$\co$, where $A \in \cc_{\co}(\tz_{J, y, \d})$. This is a
generalization of \cite[5.3]{L2}.

\subsection*{5.11} Now let us recall another definition of parabolic
character sheaves in \cite{L3}.

Let $w \in W$. The Borel group $B$ acts on $G_{I, w} \times G_{I,
y}$ by $b (g, h)=(g b \i, b h)$. The quotient $G_{I, w} \times^B
G_{I, y}$ exists. The action of $B$ on $G \times (G_{I, w} \times
G_{I, y})$ defined by $b \cdot (g, g_1, g_2)=(g b \i, b g_1, g_2
g_0 b \i g_0 \i)$ induces an action of $B$ on $G \times (G_{I, x}
\times^B G_{I, y})$. The quotient $G \times_B (G_{I, x} \times^B
G_{I, y})$ exists.

With the notation of \cite[4.2]{L2}, set
$$Y_{(w, y)}=\{(B, B', g) \in \cb \times \cb \times G^1
\mid \po(B, B')=w, \po(B', {}^g B)=y\}.$$ The morphism $G \times
(G_{I, w} \times G_{I, y}) \rightarrow Y_{(w, y)}$, $(g, g_1, g_2)
\mapsto ({}^g B, {}^{g g_1} B, g g_1 g_2 g_0 g \i)$ induces an
isomorphism $G \times_B (G_{I, w} \times^B G_{I, y}) \cong Y_{(w,
y)}$.

Let $\xi \in \hat{X}$ with $w y \d(\xi)=\xi$. Then $\tilde
\cl_{\xi, w}=\bar{\mathbb Q}_l \odot (\cl_{y \d(\xi), w} \odot
\cl_{\d(\xi), y})$ is a (tame) local system on $G \times_B (G_{I,
w} \times^B G_{I, y})$. (The local system $\tilde \cl$ on $Y_{(w,
y)}$ defined in \cite[4.2]{L2} is of the form $\tilde \cl_{\xi,
w}$ for some $\xi \in \hat{X}$ with $w y \d(\xi)=\xi$ via the
isomorphism $G \times_B (G_{I, w} \times^B G_{I, y}) \cong Y_{(w,
y)}$.)

Let $\rho_{(w, y)}: G \times_B (G_x \times^B G_{y, \d})
\rightarrow \tz_{J, y, \d}$ be the morphism defined by $(g, g_1,
g_2) \mapsto (^g P_J, {}^{g g_1} P_{J'}, g g_1 U_{P_{J'}} g_2 g_0
U_P g \i)$. Let $\cc_{J, y, \d}$ be the set of (isomorphism
classes of) simple perverse sheaves $A$ on $\tz_{J, y, \d}$ such
that $A \dashv (\rho_{(w, y)})_! \tilde \cl_{\xi, w}$ for some $w
\in W$ and $\xi \in \hat{X}$ with $w y \d(\xi)=\xi$. This is the
definition of parabolic character sheaves on $\tz_{J, y, \d}$ in
\cite[11.3]{L3}.

\begin{th17} We have that $\cc(\tz_{J, y, \d})=\cc_{J. y, \d}$.
\end{th17}

\begin{rmk} Thus we obtained a new proof of the fact that the two
definitions of parabolic character sheaves in \cite{L3} coincide,
which was first proved by Lusztig in \cite[11.15 \& 11.18]{L3}.
\end{rmk}

Proof. We identify $G_{I, w} \times^B G_{I, y}$ with $B \times_B
(G_{I, w} \times^B G_{I, y})$. Then
\begin{align*} \rho_{(w, y)} (G_{I, w} \times^B G_{I, y}) &=(P_J,
{}^{B \dot w} P_{J'}, B \dot w B \dot y B g_0) \\ &=(P_J, {}^{B
\dot w \dot y \dot y \i} P_{J'}, B \dot w \dot y U_{^{\dot y \i} P_{J'}} B g_0) \\
&=\widetilde{[J, w y, 1]}_{y, \d}.
\end{align*}

Let $\rho'_{(w, y)}: G_{I, w} \times^B G_{I, y} \rightarrow
\widetilde{[J, w y, 1]}_{y, \d}$ be the restriction of $\rho_{(w,
y)}$. Then $\rho'_{(w, y)}$ is an affine space bundle map and
$(\rho'_{(w, y)})_! (\cl_{w, y \d(\xi)} \odot \cl_{y,
\d(\xi)})=\cl_{\xi, w y, 1}$ for $\xi \in \hat{X}$. The morphism
$(id, \rho'_{(w, y)}): G \times (G_{I, w} \times^B G_{I, y})
\rightarrow G \times \widetilde{[J, w y, 1]}_{y, \d}$ induces a
morphism
$$\tilde \rho_{(w, y)}: G \times_B (G_{I, w} \times^B G_{I, y})
\rightarrow G \times_B \widetilde{[J, w y, 1]}_{y, \d}.$$

We have that $\rho_{(w, y)}=\rho \mid_{G \times \widetilde{[J, w
y, 1]}_{y, \d}} \circ \hbox{ } \tilde \rho_{(w, y)}$. Hence for
$\xi \in \hat{X}$ with $w y \d(\xi)=\xi$, we have that
$$(\rho_{(w, y)})_!  \tilde \cl_{\xi, w}=(\rho \mid_{G \times \widetilde{[J, w
y, 1]}_{y, \d}})_! (\tilde \rho_{(w, y)})_! \tilde \cl_{\xi,
w}=(\rho \mid_{G \times \widetilde{[J, w y, 1]}_{y, \d}})_!
(\bar{\mathbb Q}_l \odot \cl_{\xi, w y, 1})[\cdot].$$

Now let $A$ be a simple perverse sheaf on $\tz_{J, y, \d}$. If $A$
satisfies 5.9(iii), then $A \in \cc_{J, y, \d}$. If $A \in \cc_{J,
y, \d}$, then $A$ satisfies 5.9(ii). Therefore, $\cc(\tz_{J, y,
\d})=\cc_{J. y, \d}$. \qed

\subsection*{Note added} After completing this paper, I learned that
the main result of this paper has also been obtained by
T.A.Springer (unpublished).

\bibliographystyle{amsalpha}

\end{document}